\documentclass[psamsfonts]{article}

    \usepackage{latexsym}
    \usepackage{amsfonts}
    \usepackage{amsmath}
    \usepackage{amssymb}
    \usepackage{amscd}
    \usepackage{amstext}
    \usepackage{amsthm}
\usepackage{showkeys}

    \title{The Rokhlin property and the tracial topological rank
\thanks{Research partially supported by NSF grants DMS 0097903
and JSPS Grant for Scientific Research of Japan No.14540217.
         AMS 2000 Subject Classification Numbers:
                         Primary 46L55 and
                                46L35.
                        Key words:
                        the Rokhlin property, tracial topological rank,
simple $C^*$-algebras
                                                      \protect\\}}

    \author{Huaxin Lin\\
    Department of Mathematics\\
    University of Oregon\\
    Eugene, Oregon 97403-1222\\
    Hiroyuki Osaka\\
    Department of Mathematical Sciences\\
    Ritsumeikan University\\
    Kusatsu, Shiga, 525-8577, Japan
    }
    \date{}

\begin{document}
    \maketitle


    \newcommand{\CA}{$C^*$-algebra}
    \newcommand{\SCA}{$C^*$-subalgebra}
\newcommand{\aue}{approximate unitary equivalence}
    \newcommand{\ayue}{approximately unitarily equivalent}
    \newcommand{\mops}{mutually orthogonal projections}
    \newcommand{\hm}{homomorphism}
    \newcommand{\pisca}{purely infinite simple \CA}
    \newcommand{\andeqn}{\,\,\,\,\,\, {\rm and} \,\,\,\,\,\,}
    \newcommand{\QED}{\rule{1.5mm}{3mm}}
    \newcommand{\morp}{contractive completely
    positive linear map}
    \newcommand{\asmorp}{asymptotic morphism}
    \newcommand{\arrow}{\rightarrow}
    \newcommand{\tdsum}{\widetilde{\oplus}}
    \newcommand{\pa}{\|}  
    \newcommand{\ep}{\varepsilon}
    \newcommand{\id}{{\rm id}}
    \newcommand{\aueeps}[1]{\stackrel{#1}{\sim}}
    \newcommand{\aeps}[1]{\stackrel{#1}{\approx}}
    \newcommand{\dt}{\delta}
    \newcommand{\yu}{\fang}
    \newcommand{\ca}{{\cal C}_1}
\newcommand{\Ad}{{\rm ad}}
    \newcommand{\tr}{{\rm TR}}
    \newcommand{\Aut}{{\rm Aut}}
    \newcommand{\T}{{\rm T}}
    \newcommand{\Aff}{{\rm Aff}}
    \newcommand{\Tr}{{\rm Tr}}
    \newcommand{\Z}{{\mathbb Z}}
    \newcommand{\Q}{{\mathbb Q}}
    \newcommand{\R}{{\mathbb R}}

    \newtheorem{thm}{Theorem}[section]
    \newtheorem{Lem}[thm]{Lemma}
    \newtheorem{Prop}[thm]{Proposition}
    \newtheorem{Dfn}[thm]{Definition}
    \newtheorem{Cor}[thm]{Corollary}
    \newtheorem{Ex}[thm]{Example}
    \newtheorem{Pro}[thm]{Problem}
    \newtheorem{Remark}[thm]{Remark}
    \newtheorem{NN}[thm]{}
    \renewcommand{\theequation}{e\,\arabic{section}.\arabic{equation}}
    \newcommand{\rforal}{{\rm\,\,\,for\,\,\,all\,\,\,}}

    \newcommand{\Ik}{ {\cal I}^{(k)}}
    \newcommand{\Iz}{{\cal I}^{(0)}}
    \newcommand{\Ii}{{\cal I}^{(1)}}
    \newcommand{\Ip}{{\cal I}^{(2)}}

\begin{abstract}
Let $A$ be a unital separable simple \CA\, with $\tr(A)\le 1$ and
$\alpha$ be an automorphism.  We show that if $\alpha$ satisfies
the tracially cyclic Rokhlin property then
$\tr(A\rtimes_{\alpha}\Z)\le 1.$ We also show that whenever $A$
has a unique tracial state and $\alpha^m$ is uniformly outer for
each $m (\not= 0)$ and $\alpha^r$ is approximately inner for some $r>0,$
$\alpha$ satisfies the tracial cyclic Rokhlin property. By
applying the classification theory of nuclear \CA s, we use the
above  result to prove a conjecture of Kishimoto: if $A$ is a
unital simple $A{\mathbb T}$-algebra of real rank zero
and $\alpha\in \Aut(A)$ which is approximately
inner and if $\alpha$ satisfies some Rokhlin property, then the
crossed product $A\rtimes_{\alpha}\Z$ is again an $A{\mathbb T}$
-algebra of real rank zero.
As a by-product, we find that one can
construct a large class of simple \CA s with tracial rank one (and
zero) from crossed products.
\end{abstract}

\section{Introduction}

The Rokhlin property in ergodic theory was adopted to the context
of von Neumann algebras by A. Connes (\cite{C}). It was adopted by
Herman and Ocneanu (\cite{HO}) for UHF-algebras. M. R\o rdam
(\cite{R}) and A. Kishimoto(\cite{Ks1.5}) introduced the Rokhlin property
to a much more general context of \CA s (see also 
\cite{Iz}). A. Kishimoto had been studying
automorphisms on UHF-algebras and more generally, on simple
$A{\mathbb T}$-algebras that satisfy a Rokhlin property
 (\cite{Ks1}, \cite{Ks2}).
More recently,
 N. C. Phillips studied finite group
actions which satisfy certain type of Rokhlin property on some
simple \CA s (\cite{P}). 

A conjecture  of A. Kishimoto can be formulated as follows: Let $A$
be a unital simple $A{\mathbb T}$-algebra of real rank zero and
$\alpha$ be an approximately inner automorphism. Suppose that
$\alpha$ is ``sufficiently outer", then the crossed product of the
$A{\mathbb T}$-algebra by $\alpha,$ $A\rtimes_{\alpha}\Z$ is again
a unital simple $A{\mathbb T}$-algebra. In particular, he studied
the case that $A$ has a unique tracial state.

A. Kishimoto proposed that the appropriate notion of the outerness is
the Rokhlin property (\cite{Ks2}). He also introduced the notion
of uniformly outer (\cite{Ks1.5}). In \cite{Ks2}, he showed that
if $A$ is a unital simple $A{\mathbb T}$-algebra of real rank zero
with a unique tracial state and $\alpha\in \Aut(A)$ is approximately
inner, then $\alpha$ has the Rokhlin property if and only if
$\alpha^m$ is uniformly outer for all $m\not=0.$ He also showed
that the Rokhlin property, in this situation, is  equivalent to
say that $A\rtimes_{\alpha}\Z$ has real rank zero, and it is
equivalent to say that $A\rtimes_{\alpha}\Z$ has a unique tracial
state. A. Kishimoto showed (\cite{Ks1})  that, if in addition, $A$ is
a UHF-algebra then $A\rtimes_{\alpha}\Z$ is in fact a unital
simple $A{\mathbb T}$-algebra. He also showed in \cite{Ks3} that
the conjecture is true for the case that $A$ is assumed to have a
unique tracial state, both $K_i(A)$ are finitely generated and
$K_1(A)\not\cong \Z,$ and, in addition, that $\alpha\in {\rm
Hinn}(A).$ Among other things, we prove in this paper
 the Kishimoto conjecture for all cases that $A$ has a unique
tracial state. If the term of ``sufficiently outer" is interpreted
as ``tracial Rokhlin" property, then Kishimoto's conjecture holds:
Let $A$ be a unital simple $A{\mathbb T}$-algebra and $\alpha$ be
an approximate inner automorphism. Suppose that $\alpha$ has the
tracial Rokhlin property, then $A\rtimes_{\alpha}{\mathbb T}$ is
again a unital simple $A{\mathbb  T}$-algebra.

We take the advantage of the development in Elliott's program
of  the classification of nuclear \CA s (see \cite{Ell1} and
\cite{EG}, for example). In particular, we use the classification result in
\cite{Lnmsri}, where unital separable simple \CA s satisfying the
Universal Coefficient Theorem and with tracial topological rank zero are
classified by their $K$-theory. Adopting a N. C. Phillips's
observation, we note that if $A$ is a unital simple \CA\, with
$\tr(A)\le 1$ and $\alpha\in \Aut(A)$ satisfies a so-called tracial
cyclic Rokhlin property, then $\tr(A\rtimes_{\alpha}\Z)\le 1$ so
that the classification result in \cite{Lnmsri} and \cite{Lntai}
can be applied. Using Kishimoto's techniques, we show that if
$\alpha^r$ is approximately inner (for some integer $r>0$), the
tracial Rokhlin property introduced in (\cite{OP}) implies the
tracial cyclic Rokhlin property. Using a result in \cite{OP}, we
actually show a more general result (see Theorem \ref{IIT2}).

The assumption that $\alpha$ is approximately inner is to insure
that the crossed products remain finite (see the introduction of
\cite{Ks2}). We relax this restriction slightly by only requiring
that $\alpha^r$ is approximately inner for some integer $r>0.$ It
turns out that in a number of cases, while there are automorphisms
which are not approximate inner, all (outer) automorphisms
$\alpha$ have this property, i.e., for some integer $r>0,$
$\alpha^r$ are approximately inner (see Theorem \ref{IVTainn}).
We show that our results also cover many cases in which $A$ may have
arbitrary tracial space (see Corollary \ref{IVCtr=1}).

It is shown by G. Gong (\cite{G}) that a unital simple AH-algebra
with very slow dimension growth has tracial topological rank one or zero.
Moreover, Elliott, Gong and Li (\cite{EGL})  show that the class
of unital simple AH-algebras with very slow dimension growth can
be classified by their $K$-theoretical data. An improvement of
this classification has been made so that unital simple nuclear
\CA s with tracial topological rank no more than one which satisfy the
Universal Coefficient Theorem can also be classified by their
$K$-theoretical data (\cite{Lntai}). However, until now, all
interesting examples of unital simple nuclear \CA s that have
tracial topological rank one are those AH-algebras with very slow dimension
growth (and those of similar inductive limit construction).
Theorem \ref{IT1} also provides ways
to construct unital simple \CA s with tracial topological rank one by crossed
products (see Corollary \ref{IVCtr=1} and Example \ref{IVE}).
It also creates  the
opportunity to apply the classification results in \cite{Lntai}.

{\bf Acknowledgement} Most of this work was done when the first
author was visiting Ritsmeikan University. He acknowledges the
support of JSPS Grant for Scientific Research of Japan.
During the work he was also supported by
National Science Foundation of USA.
The second author also acknowledges the
support of JSPS Grant for Scientific Research of Japan.
Both authors would like to thank Masaru Nagisa for a fruitful discussion.

\section{The Rokhlin properties}

The following conventions will be used in this paper.
Let $A$ be a unital \CA~

\begin{itemize}
\item[(i)]
We denote by $\Aut(A)$ the set of all automorphisms on $A$ and by
$\T(A)$ the tracial state space of $A$.
\item[(ii)]
Two projections $p, q \in A$ are said to be equivalent if
they are Murray-von Neumann equivalent. That is, there exits
a partial isometry $w \in A$ such that $w^*w = p$ and $ww^* = q$.
Then We write $p \sim q$.
\item[(iii)]
Let  ${\mathcal F}$ and ${\mathcal S}$ be subsets of $A$ and
$\ep > 0$.
We write $x \in_\ep {\mathcal S}$ if there exists $y \in {\mathcal S}$
such that $\|x - y\| < \ep$, and
write ${\mathcal F} \subset_\ep {\mathcal S}$ if $x \in_\ep {\mathcal S}$
for all $x \in {\mathcal F}$.
\item[(iv)]
Let $a$ and $b$ be two positive elements in $A$. We write $[a]
\leq [b]$ if there exists an element $x \in A$ such that $a =
x^*x$ and $xx^* \in \overline{bAb}.$
If $ab = ba = 0$, then we write $[a + b] = [a] + [b]$. Let $p$ be
a projection and $b$ non-zero positive element in $A$. Note that
$[p] \leq [b]$ implies that  $p$ is Murray-von Neumann equivalent
to a projection in the hereditary \CA~ $\overline{bAb}$.
\item[(v)]
We denote by ${\mathcal I}^{(0)}$ the class of all finite dimensional
C*-algebras, and by ${\mathcal I}^{(k)}$
the class of all C*-algebras with the form $pM_n(C(X))p$, where
$X$ is a finite CW complex with dimension $k$ and $p \in M_n(C(X))$
is a projection.
\end{itemize}

\vskip 3mm

We recall the definition of tracial topological rank of
C*-algebras.

\begin{Dfn}\label{ID0}{\rm \cite[Theorem 6.13]{Lntr0}}
{\rm
Let $A$ be a unital simple \CA~ and $k \in {\mathbb N}$. Then
$A$ is said to  have
{\it tracial topological rank no more than $k$} if and only if
for any finite set ${\mathcal F} \subset A$, and $\ep > 0$ and
any non-zero positive element $a \in A$, there exists a
\SCA\, $B \subset A$ with $B \in {\mathcal I}^{(k)}$ and
${\rm id}_B = p$ such that
\begin{itemize}
\item[$(1)$]
$\|[x, p]\| < \ep$ for all $x \in {\mathcal F}$,
\item[$(2)$] $pxp \in_\ep B$ for all $x \in {\mathcal F}$,
\item[$(3)$]
$[1 - p] \leq [a]$.
\end{itemize}
We write $\tr(A) \leq k$.
}
\end{Dfn}

\vskip 3mm

\begin{Remark}\label{IR0}{\rm \cite[Corollary 6.15]{Lntr0}}
{\rm Let $A$ be a simple unital C*-algebra with stable rank one
which satisfies the Fundamental Comparison Property. 
Then $\tr(A)
\leq k$ if and only if for any finite set ${\mathcal F} \subset
A$, $\ep > 0$, and any non-zero positive element $a \in A$, there
exists a \SCA\, $B \subset A$ with $B \in {\mathcal
I}^{(k)}$ and ${\rm id}_B = p$ such that
\begin{itemize}
\item[$(1)$]
$\|[x, p]\| < \ep$ for all $x \in {\mathcal F}$,
\item[$(2)$] $pxp \in_\ep B$ for all $x \in {\mathcal F}$,
\item[$(3)$]
$\tau(1 - p) < \ep$ for all $\tau \in \T(A)$.
\end{itemize}
Recall that $A$ is said to have
{\it the Fundamental Comparison Property} if $p, q \in A$
are two projections with $\tau(p) < \tau(q)$ for all
$\tau \in \T(A)$, then $p$ is equivalent to a subprojection of
$q$.
}

\end{Remark}

\vskip 3mm

The following is defined in \cite[Definition 2.1]{OP}.

\begin{Dfn}\label{ID1}
{\rm
Let $A$ be a simple unital \CA~ and
let $\alpha \in \Aut(A)$.
We say $\alpha$ has the {\it tracial Rokhlin property}
if for every finite set $F \subset A$, every $\ep > 0$,
every $n \in {\mathbb N}$, and every nonzero positive element $x \in A$,
there are mutually orthogonal projections $e_0, e_1, \dots, e_n \in A$
such that:
\begin{itemize}
\item[$(1)$]
$\| \alpha (e_j) - e_{j + 1} \| < \ep$ for $0 \leq j \leq n - 1$.
\item[$(2)$]
$\| e_j a - a e_j \| < \ep$ for $0 \leq j \leq n$ and all $a \in F$.
\item[$(3)$]
With $e = \sum_{j = 0}^{n} e_j$, $[1 - e] \leq [x]$.
\end{itemize}
}
\end{Dfn}

\vskip 3mm

We define a slightly stronger version of the tracial Rokhlin
property similar to the approximately Rokhlin property in
\cite[Definition 4.2]{Ks1}.

\vskip 3mm

\begin{Dfn}\label{ID2}
{\rm
Let $A$ be a simple unital
 \CA~ and
let $\alpha \in \Aut(A)$.
We say $\alpha$ has the {\it tracial cyclic Rokhlin property}
if for every finite set $F \subset A$, every $\ep > 0$,
every $n \in {\mathbb N}$, and every nonzero positive element $x \in A$,
there are mutually orthogonal projections $e_0, e_1, \dots, e_n \in A$
such that
\begin{itemize}
\item[$(1)$]
$\| \alpha (e_j) - e_{j + 1} \| < \ep$ for $0 \leq j \leq n$,
where $e_{n + 1} = e_0$.
\item[$(2)$]
$\| e_j a - a e_j \| < \ep$ for $0 \leq j \leq n$ and all $a \in F$.
\item[$(3)$]
With $e = \sum_{j = 0}^{n} e_j$, $[1 - e] \leq [x]$.
\end{itemize}
}
\end{Dfn}

\vskip 3mm

\begin{Remark}\label{IR1}
{\rm
\begin{itemize}
\item[(i)]
The only difference between the tracial Rokhlin property and the
tracial cyclic Rokhlin property is that in condition $(1)$ we require
that  $\|\alpha(e_n) - e_0\| < \ep.$
\item[(ii)]
If $A$ has real rank zero, stable rank one and has  weakly
unperforated $K_0(A)$ (or if $A$ has SP-property, stable rank one,
 and the Fundamental Comparison Property),
then condition $(3)$ in both Rokhlin
property can be replaced by the following condition $(3)'$
using the standard argument:

$(3)'$
With $e = \sum_{j=0}^ne_j$, we have
$\tau(1 - e) < \ep$ for all $\tau \in \T(A)$.
\end{itemize}

\item[(iii)] If $A$ is a simple unital \CA~
with real rank zero, stable rank one, and has weakly
unperforated $K_0(A)$, the Rokhlin property in the sense
of Kishimoto (\cite{Ks1}) implies the tracial Rokhlin
property (\cite{OP}).

\vskip 3mm

Recall that a C*-algebra $A$ is said to have {\it SP-property}
if any non-zero hereditary C*-subalgebra of $A$ has
a non-zero projection. 
}
\end{Remark}

\vskip 3mm

Obviously the tracial cyclic Rokhlin property implies the tracial
Rokhlin property. The converse is also true in many cases. We
will discuss it in the next section.

Before stating the characterization of the tracial Rokhlin
property we cite the following notion introduced by Kishimoto
(\cite{Ks1.5}).

\begin{Dfn}\label{ID3}
{\rm
Let $A$ be a unital \CA~ and $\alpha \in \Aut(A)$.
We say $\alpha$ is uniformly outer if for any $a \in A$,
any projection $
p \in A$, and any $\ep > 0$, there are
finite number of projections $p_1, \dots, p_n$in $A$
such that $\sum_ip_i = p$ and $\|p_ia\alpha(p_i)\| < \ep$
for $i = 1, \dots, n$.
}
\end{Dfn}

\vskip 3mm

The following result is the tracial Rokhlin version of
Kishimoto's result in the case of simple unital $A{\mathbb T}$-algebras
with a unique trace \cite[Theorem 2.1]{Ks1}.

\begin{thm}\label{IT1}{\rm \cite{OP}}
{\rm
Let $A$ be a simple unital \CA~ with $\tr(A) = 0,$
and suppose that $A$ has a unique tracial state.
Then the following conditions are equivalent:
\begin{itemize}
\item[$(1)$]
$\alpha$ has the tracial Rokhlin property.
\item[$(2)$]
$\alpha^m$ is not weakly inner in the GNS representation
$\pi_\tau$ for any $m \not= 0$.
\item[$(3)$]
$A \rtimes_\alpha{\mathbb Z}$ has real rank zero.
\item[$(4)$] $A \rtimes_\alpha{\mathbb Z}$ has a unique trace.
\end{itemize}
}
\end{thm}

\vskip 3mm

Note that the uniformly outerness implies that $\alpha$ is
 not weakly inner in the GNS representation $\pi_\tau$ by
 an $\alpha$-invariant tracial state $\tau$ on $A$
 by \cite[Lemma 4.4]{Ks1.5}.

\begin{Remark}\label{IR2}
{\rm When $A$ is  a simple unital \CA~ with tracial topological
rank zero, if $\alpha \in \Aut(A)$ has the tracial Rokhlin
property, it is proved in \cite{OP} that the crossed product $A
\rtimes_\alpha{\mathbb Z}$ has real rank zero, stable rank one,
and the order on projections over $A \rtimes_\alpha{\mathbb Z}$ is
determined by traces. But it is not known that the crossed product
$A \rtimes_\alpha{\mathbb Z}$ has tracial topological rank zero.
However, if $\alpha$ has the tracial cyclic Rokhlin property, then we have
the following result based on an observation of N. C. Phillips
(\cite{P}). }
\end{Remark}

\vskip 3mm

\begin{thm}\label{IT2}
{\rm Let $A$ be a simple unital \CA~ with $\tr(A) \leq 1.$ Suppose
that  $\alpha \in \Aut(A)$ has the tracial cyclic Rokhlin property.
Then $\tr(A \rtimes_\alpha{\mathbb Z}) \leq 1$.

In particular, if $A$ has $\tr(A) = 0$, then
$\tr(A \rtimes_\alpha{\mathbb Z}) = 0$.
}
\end{thm}

\begin{proof}
We first note that, by \cite{Ks0},
$A\rtimes_{\alpha}{\mathbb Z}$ is a simple \CA.

Let $\ep > 0$, $n \in {\mathbb N}$, and $F \subset A \rtimes_\alpha{\mathbb Z}$
be a finite set.
To simplify notation, without loss of generality, we may assume that
$$
F = \{a_i\}_{i=1}^m \cup \{u\},
$$
where $a_i\in A$ and $\|a_i\|\le 1$ ($i=1,2,...,m$ ) and $u$ is a
unitary which implements $\alpha$. Fix $b\in (A
\rtimes_\alpha{\mathbb Z})_+\backslash \{0\}.$

Since $A$ has SP-property and $\alpha$ is outer, $A
\rtimes_\alpha{\mathbb Z}$ also has  SP-property \cite[Theorem
4]{JO}. In particular, there is a  non-zero  projection $r \in
\overline{b(A \rtimes_\alpha{\mathbb Z})b}.$ Let $r_0\in A$ be a
nonzero projection. Since $A\rtimes_{\alpha}\Z$ is simple, by 1.8 of
\cite{Cu}, it is easy to find  a non-zero projection $r' \in
r_0Ar_0$ such that $r'$ is equivalent to a subprojection of $r$
(see, for example, \cite[Theorem 4]{JO}).
Hence there are projections $r_1, r_2 \in A$ such that $r_1r_2 =
0$ and $r_1 + r_2$ is equivalent to a subprojection of $r$
(see, for example, 3.5.7 of \cite{Lnbk}).

Since $\alpha$ has the tracial cyclic Rokhlin property, for any
$\delta > 0$ with $\delta <
{{\ep}\over{5}}$
there exist projections $e_1, e_2$ such that
\begin{itemize}
\item[(1)]
$\|\alpha(e_i) - e_{i+1}\| < \delta$ for $1 \leq i \leq 2$
($e_{3} = e_1$)
\item[(2)]
$\|[e_i, a_k]\| < \delta$ for $1 \leq k \leq m$.
\item[(3)]
$[1 - e_1 - e_2] \leq [r_1]$.
\end{itemize}

Set $p = e_1 + e_2$. From (1) above, one estimates that
$$
\begin{array}{ll}
\|up - pu\| &= \|\sum_{i=1}^2ue_i - \sum_{i=1}^2e_{i+1}u\|\\
&= \sum_{i=1}^2\|ue_i - e_{i+1}u\| < 2\delta.
\end{array}
$$
Hence, together with (2) above, we obtain
\begin{itemize}
\item[(4)]
$\|[p, a]\| < 2\delta\rforal a\in F.$
\end{itemize}
There is a unitary $v \in A \rtimes_\alpha{\mathbb Z}$ such that
$\|v - 1\| < \delta$ and $vu^*e_iuv^* = e_{i+1}$ for $1 \leq i
\leq 2$. Set $w = vu^*$, and consider the C*-algebra $D$ generated by
$e_1Ae_1$ and $e_2we_1$. Then $D$ is isomorphic to $e_1Ae_1
\otimes M_2({\mathbb C})$. Note that $pw=e_1w+e_2w=we_2+we_1=wp.$
Moreover, $pwp\in D.$ Since $\|pup-pwp\|<\dt,$ one has that
$pup\in_{\dt} D.$ By (2) again, we have
$$
\|pa_jp-(e_1a_je_1+e_2a_je_2)\|<2\dt\,\,\,j=1,2,...,m.
$$
It follows that
$pFp \subset_{2\delta}D$.

Since $A$ is a simple \CA~ with SP-property, there exists a
nonzero projection $r_3 \in e_1Ae_1$ such that $r_3$ is equivalent
to a subprojection of  $r_2.$ Since $\tr(e_1Ae_1)\le 1$
($\tr(e_1Ae_1)=0$ if $\tr(A)=0$), $\tr(D)\le 1$ ($\tr(D)=0$ if
$\tr(A)=0$) by \cite[Theorem 5.3]{Lntr0}.  So there exists a
\SCA\, $B \in {\mathcal I}^{(k)}$ \ ($k = 1$ or $k = 0$)
 and projection $e = 1_B$
such that
\begin{itemize}
\item[(5)]
$\|[pap, e]\| < \delta < \ep\rforal a\in F,$
\item[(6)]
$G \subset_\delta B\andeqn $
\item[(7)]
$[p - e] \leq [r_3]$.
\end{itemize}

{}From (3), (4), (5), and (7) above
we estimate that
\begin{itemize}
\item[(8)]
For any $f \in F$
$$
\begin{array}{ll}
\|ef - fe\| &=\|e(pf - fp) + efp - pfe + (pf - fp)e\|\\
&\leq \|pf - fp\| + \|epfp - pfpe\| + \|pf - fp\|\\
&< 5\delta < \ep \andeqn
\end{array}
$$
\item[(9)]
$$
\begin{array}{ll}
[1 - e] &= [1 - p + p - e],\\
&= [
1 - p] + [p - e]\\
&\leq [r_1] + [r_3]  \leq [r_1] + [r_2] \leq [r] \leq [b].
\end{array}
$$
\end{itemize}

{}From (6) and $pFp \subset_{2\delta} D$, we have

(10) $pFp \subset_{4\delta} D$.

Hence from the estimates (8), (9), (10) we conclude that
$A \rtimes_\alpha{\mathbb Z}$ has tracial topological
less than or equal to $1$.

In the case of $\tr(A) = 0$ $B$ can be chosen to be finite dimensional.
Hence, in that case, $\tr(A \rtimes_\alpha{\mathbb Z}) = 0$.
\end{proof}


\section{Approximately inner automorphisms}

\begin{Lem}\label{IIL1}
{\rm
Let $A$ be a unital separable \CA\,
and $\alpha: A\to A$ be an approximate inner automorphism. Suppose
that $\{p_j\}$ is a central sequence of projections. Then there
exists a central sequence of partial isometries $\{w_j\}$ such
that $w_j^*w_j = p_j$ and $w_jw_j^* = \alpha(p_j),$ $j=1, 2,....$
}
\end{Lem}

\begin{proof}
Fix a finite subset ${\cal F}\subset A$ which is in the unit
ball of $A.$  Let $\ep>0.$ Choose a unitary $v\in U(A)$ such that
$$
\|\alpha(a)-v^*av\|<\ep/8\rforal a\in {\cal F}.
$$
Since $\alpha$ is an automorphism, $\alpha(p_j)$ is also a central
sequence of projections.
Choose a sufficiently large $j$ so that
$$
\|p_ja-ap_j\|<\ep/8\rforal a\in {\cal F}\andeqn
\|\alpha(p_j)v-v\alpha(p_j)\|<\ep/8.
$$
Since $\alpha$ is approximately inner, we obtain another unitary
$z\in U(A)$ such that
$$
\|z^*p_jz-\alpha(p_j)\|<\ep/8\andeqn
\|z^*az-\alpha(a)\|<\ep/8\rforal a\in {\cal F}.
$$
It follows that
$$
\|(vz^*)p_j(zv^*)-\alpha(p_j)\|\le
\|vz^*p_jzv^*-v\alpha(p_j)v^*\|+\|v\alpha(p_j)v^*-\alpha(p_j)\|<\ep/4.
$$
Let $x_j=vz^*p_j.$ Then $x_j^*x_j=p_j$ and
$$
\|x_jx_j^*-\alpha(p_j)\|<\ep/4.
$$
{}From the above we also have
 $$
\|vz^*azv^*-a\| < \ep/4\andeqn
\|vz^*a-avz^*\| < \ep/4\rforal a\in {\cal F}.
$$

 On the other hand, for any $a \in {\cal F},$
$$
\|x_ja-ax_j\| \le
\|vz^*p_ja-vz^*ap_j\| + \|vz^*ap_j-avz^*p_j\| < \ep/8+\ep/4 = 3\ep/8.
$$
There is a unitary $u \in U(A)$ such that
$\|u-1\| < \ep/4$ such that

$$
u(x_jx_j^*)u^* = \alpha(p_j).
$$

Define $w_j = ux_j.$ Then $w_jw_j^*=p_j$ and $w_jw_j^*=\alpha(p_j).$
Moreover we have that

$$
\|w_ja-aw_j\|<\ep\rforal a\in {\cal F}.
$$

Since ${\cal F}$ is arbitrary, the lemma follows.
\end{proof}

\begin{Lem}\label{IIL2}
{\rm
Let $A$ be a unital separable \CA\,
and $\alpha\in \Aut(A)$ for which
$\alpha^r$ is approximately inner for some integer $r\ge 1.$ Let
$m\in {\mathbb N},$  $m_0\ge m$ be the smallest integer such that
$m_0 = 0 \ {\rm mod}\, r$ and $l = m+(r-1)(m_0+1).$

Suppose that $\{e_i^{(n)}\},$ $i = 0,1,...,l,$ $n = 1, 2,...,$ are
$l+1$ sequences of projections in $A$ satisfying the following:

$$
\|\alpha(e_i^{(n)})-e_{i+1}^{(n)}\|<\dt_n, \,
\lim_{n\to\infty}\dt_n=0,
$$

$$
e_i^{(n)}e_j^{(n)}=0,\, {\rm if}\,\,\, i\not=j,\, e_i^{(n)}\sim
e_j^{(n)}\,\,\,{\rm in} \,\,\, A,
$$

and for each $i,$ $\{e_i^{(n)}\}$ is a central sequence.

Then for each $i = 0, 1, 2, \dots, m,$
there is a central sequence of partial
isometries $\{w_i^{(n)}\}$ such that

$$
(w_i^{(n)})^*w_
i^{(n)} = p_{i}^{(n)}\andeqn
w_i^{(n)}(w_i^{(n)})^* = p_{i+1}^{(n)},\,i=0,1,...,m - 1,
$$

where $p_i^{(n)} = \sum_{j=0}^{r-1}e_{i+j(m_0+1)}^{(n)}.$
Moreover, for each $i$,

$$
\lim_{n \rightarrow \infty}\|\alpha(p_i^{(n)}) - p_{i+1}^{(n)}\| = 0.
$$
}
\end{Lem}

\begin{proof}
Since $\alpha^r$ is approximately inner, by applying Lemma \ref{IIL1},
for each $i,j = 0, 1,... l,$ one obtains central sequences of partial
isometries $\{z(i,j,n)\}$ such that

$$
z(i,j,n)^*z(i,j,n) = e_{i}^{(n)}\andeqn
z(i,j,n)z(i,j,n)^* = \alpha^{rj}(e_i^{(n)}).\,
$$

Note that
$$
\|\alpha^{rj}(e_i^{(n)})-e_{i+rj}^{(n)}\|<rj\dt_n.
$$
There is a unitary $u(i,j,n) \in U(A),$ for each $i$ and $j,$ such
that

$$
\|u(i,j,n)-1\| < 2(rl)\dt_n
\andeqn
u(i,j,n)^*\alpha^{rj}(e_{i}^{(n)})u(i,j,n)
= e_{i+rj}^{(n)}
$$

Since $\lim_{n\to\infty} \dt_n=0,$ for each $i$ and $j,$ $\{u(i,j,n)\}$ is central.
Therefore, to simplify notation, we may assume that

$$
z(i,j,n)^*z(i,j,n) = e_{i}^{(n)}\andeqn z(i,j,n)z(i,j,n)^* =
e_{i+rj}^{(n)}.
$$

Define

$$
p_i^{(n)} = \sum_{j=0}^{r - 1}e_{i+j(m_0+1)}^{(n)}.
$$

Then for each $i,$
one checks easily that there are central sequences of partial isometries
$\{w(i,n)\}$ such that

$$
w(i,n)^*w(i,n) = p_i^{(n)}\andeqn w(i,n)w(i,n)^* = p_{i+1}^{(n)},\,i=0,1,...,m.
$$

For example, (with $m_0=kr$), one defines
$$
w(0,n) = z(1,k,n)^* + z(1+(m_0+1),k,n)^* + \cdots +
z(1+(r-2)(m_0+1),k,n)^* + z(0,(r-1)k+1,n).
$$
Then (with $m_0 = kr$)
$$
\begin{array}{ll}
&w(0,n)w(0,n)^* \\
&= z(1,k,n)^*z(1,k,n)+
z(1+(m_0+1),k,n)^*z(1+(m_0+1),k,n)
+ \cdots \\
&+ z(1+(r-2)(m_0+1),k,n)^*z(1+(r-2)(m_0+1),k,n)\\
&+ z(0,(r-1)k+1,n)z(0,(r-1)k+1,n)^*\\
&=
e_1^{(n)} + e_{1+(m_0+1)}^{(n)} + \cdots +
e_{1+(r-2)(m_0+1)}^{(n)} +
e_{1+(r-1)(m_0+1)}^{(n)}\\
&= p_1^{(n)}
\end{array}
$$
(note that  $e_{((r-1)k+1)r}^{(n)} = e_{(r-1)kr+r}^{(n)} =
e_{(r-1)m_0+r}^{(n)} = e_{1+(r-1)(m_0+1)}^{(n)}$)
and
$$
\begin{array}{ll}
&w(0,n)^*w(0,n) \\
&= z(1,k,n)z(1,k,n)^* +
z(1+(m_0+1),k,n)z(1+(m_0+1),k,n)^*
+ \cdots \\
&+ z(1+(r-2)(m_0+1),k,n)z(1+(r-2)(m_0+1),k,n)^*\\
&+ z(0,(r-1)k+1,n)^*z(0,(r-1)k+1,n)\\
&= e_{1+kr}^{(n)}+e_{1+(m_0+1)+kr}^{(n)}+ \cdots +
e_{1+(r-2)(m_0+1)+kr)}^{(n)}+
e_{0}^{(n)}\\
&=e_{m_0+1}^{(n)}+e_{2(m_0+1)}^{(n)}+\cdots+e_{(r-1)(m_0+1)}^{(n)}+e_0^{(n)}\\
&= p_0^{(n)},\,\,\,\andeqn \\
\end{array}
$$
Since, for each
$i$ and $j,$  $\{z(i,j,n)\}$ is central, so is $\{w(i,n)\}.$

{}From the construction we know that for each $i$

$$
\begin{array}{ll}
\|\alpha(p_i^{(n)}) - p_{i+1}^{(n)}\| &\leq
\sum_{j=0}^{r - 1}\|\alpha(e_{i+j(m_0 + 1)}) - e_{i + 1 +j(m_0 + 1)})\|\\
&\leq r\delta_n \rightarrow 0 \ (n \rightarrow \infty)
\end{array}
$$
\end{proof}

\vskip 3mm

Let $\{E_{i,j}\}$ be a system of matrix units and ${\cal K}$ be
the compact operators on $\ell^2({\mathbb Z})$ where we identify
$E_{i,i}$ with the one-dimensional projection onto the functions
supported by $\{i\} \subset {\mathbb Z}$.
Let $S$ be the canonical shift operator on $\ell^2({\mathbb Z})$.
Define an automorphism $\sigma$ of
${\cal K}$ by $\sigma(x) = SxS^*$ for all $x \in {\cal K}$.
Then $\sigma(E_{i,j}) = E_{i + 1, j + 1}$.
For any $N \in {\mathbb N}$
let $P_N = \sum_{i = 0}^{N - 1}E_{i,i}$.

\begin{Lem}{\rm (Kishimoto, 2.1 of \cite{Ks1})}\label{LK}
{\rm For any $\eta > 0$ and $n \in {\mathbb N}$ there exist $N \in
{\mathbb N}$ and projections $e_0, e_1, \dots, e_{n - 1}$ in
${\cal K}$ such that

$$
\begin{array}{l}
\sum_{i=0}^{n - 1}e_i \leq P_N\\
\|\sigma(e_i) - e_{i + 1}\| < \eta, \ i = 0, \dots, n - 1, \ e_n = e_0\\
\frac{n\dim e_0}{N} > 1 - \eta.
\end{array}
$$
}
\end{Lem}

\begin{thm}\label{IIT1}
{\rm Let $A$ be a unital separable simple \CA\, with $\tr(A)\le 1$
and $\alpha^r \in \Aut(A)$ be an approximately inner automorphism
for some integer $r \geq 1$. Suppose that $\alpha$ has the tracial
Rokhlin property then $\alpha$ has the tracial cyclic Rokhlin
property. }
\end{thm}

\begin{proof}
Let $\ep>0.$  Let $\ep/2 > \eta > 0$ and $m \in {\mathbb N}$ be
given. Choose $N$ which satisfies the conclusion of Lemma \ref{LK}
(with this $\eta$ and $n=m$). Identify $P_N{\cal K}P_N$ with
$M_N.$ Let ${\cal G}=\{E_{i+1,i}:i=0,1,...,N-1\}$ be a set of
generators of $M_N.$ Let $e_0, e_1,...,e_{m-1}$ be as in the
conclusion of Lemma \ref{LK}.
%

 For any $\ep > 0,$
there is $\dt > 0$ depends only on $N$
such that, if

$$
\|ag-ga\| <\dt
$$
for $g \in {\cal G}$,
then

$$
\|ae_i - e_ia\| < \ep/2, \, i = 0, 1,..., n.
$$

We assume that $\dt < \eta.$  Fix a finite subset ${\cal
F}_0\subset A.$

Choose $m_0 \in {\mathbb N}$ such that $m_0 \geq m$ is the
smallest integer with $m_0 = 0 \ {\rm mod}\ r$.
Let $L = N + (r - 1)(m_0 + 1)$.

Since $\alpha$ has the tracial Rokhlin property, there exits
a sequence of projections
$\{e_i^{(k)}: i = 0, 1,...,L\}$
satisfying the following:

$$
\|\alpha(e_i^{(k)})-e_{i+1}^{(k)}\|<{\dt\over{(2^k)4N}},\,\,\,
e_i^{(k)}e_j^{(k)}=0,\,\,\,
{\rm if}\,\,\, i \not =  j,
$$

$$
\lim_{k\to\infty}\|e_i^{(k)}a-ae_i^{(k)}\| = 0
\rforal a \in A,\,i = 0, 1,..., L \andeqn
$$

$$
\tau(1-\sum_{i = 0}^{L-1}e_i^{(k)}) < \eta \rforal \tau \in \T(A),\
k = 1, 2,....
$$



By applying Lemma \ref{IIL2},
we obtain a central sequence $\{w_i^{(k)}\}$ in $A$
such that

$$
\begin{array}{l}
(w_i^{(k)})^*w_i^{(k)} = P_0^{(k)}   \andeqn\\
w_i^{(k)}(w_i^{(k)})^* = P_i^{(k)}, \ k = 0, 1,...,\
i = 0, 1, \dots, N, \\
P_i^{(k)}P_j^{(k)} = 0 \ i \not= j,\\
\|\alpha(P_i^{(k)}) - P_{i + 1}^{(k)}\| < \frac{\dt}{4L},
\ k = 0, 1,...,\ i = 0, 1, \dots, N - 1,\\
\tau(1 - \sum_{i=  0}^{N-1}P_i^{(k)}) < \eta, \ \rforal \tau \in \T(A)
\end{array}
$$
where $P_i^{(k)} = \sum_{j = 0}^{r-1}e_{i+(m_0+1)j}^{(k)}$
for $i = 0, 1, \dots, N$.

It follows that $\{ \alpha^l(w^{(k)} )\},$ $l = 0, 1,..., N$ are
all central sequences. As the same argument in Lemma \ref{IIL2}
there is a unitary $u_k \in U(A)$ with $\|u_k - 1\|<\dt/2N$ such
that ${\rm ad}\, u_k \circ \alpha(P_i^{(k)}) = P_{i+1}^{(k)},$ $i
= 0, 1,..., N-1.$ Put $\beta_k = {\rm ad}\, u_k\circ \alpha$, and
$w^{(k)} = w_0^{(k)}$. Choose a large $k,$ such that

$$
\|\beta_k^{l}(w^{(k)})a - a\beta_k^{l}(w^{(k)})\| < \dt\rforal
a\in {\cal F}_0,
$$
$l = 0, 1, ..., N.$

Now let $C_1$ and  $C_2$
be the \CA s generated by
$w^{(k)}, \beta_k^{1}(w^{(k)}), ..., \beta_k^{N - 1}(w^{(k)})$ and by
$w^{(k)}, \beta_k^{1}(w^{(k)}), ..., \beta_k^{N}(w^{(k)}),$
respectively. Note that $C_1\cong M_{N},$ $C_2\cong M_{N + 1}.$
Define a \hm\, $\Phi: C_1\to {\cal K}$ by

$$
\Phi(\beta_k^{i}(w^{(k)})) = E_{i+1,i}, \ i = 0, 1, ..., N - 1
$$

(see Lemma \ref{LK}).
Then one has $\sigma\circ \Phi|_{C_1} = \Phi\circ
\beta_k|_{C_1}$ and $\Phi(C_1) = P_{N}{\cal K}P_{N}.$
Now we apply Lemma \ref{LK} to obtain
mutually orthogonal projections $e_0, e_1,...,e_{m-1}$
in $M_N$ such that
$$
\|\sigma(e_i)-e_{i-1}\|<\eta\andeqn
{m\dim e_0\over{N}} > 1 - \eta.
$$
Let
$p_i = \Phi^{-1}(e_i),$ $i = 0, 1,..., m - 1.$
One estimates that

$$
\tau(\sum_{i = 0}^{N-1}P_i^{(k)} - \sum_{i=0}^{m-1}p_i) <
1 - \sum_{i=0}^{m-1}\frac{\dim(e_0)}{N}
= 1 - \frac{m\dim(e_0)}{N} < \eta < \frac{\ep}{2}
$$
for all $\tau \in T(A).$
So one has mutually orthogonal projections
$p_0, p_1, p_2,..., p_{m - 1}$ such that
$$
\|\beta_k(p_i) - p_{i+1}\| < \frac{\ep}{2},
 i = 0, 1, 2,..., m - 1,\, p_m = p_0.
$$
By the choice of $\dt,$ one also has
$$
\|ap_i - p_ia\| < \ep, \  i = 0, 1,..., m - 1, \ \rforal a \in
{\cal F}_0
 \andeqn
$$

$$
\tau(1-\sum_{i=0}^{m - 1}p_i) < \tau(1 - \sum_{i=0}^{N-1}P_i^{(k)})
+ \frac{\ep}{2} < \eta + \frac{\ep}{2}< \ep,
$$
for all $\tau \in \T(A).$
Since
$$
\|\beta_k - \alpha\| < \dt/2 < \ep/2,
$$
one finally has
$$
\|\alpha(p_i) - p_{i+1}\| < \ep,\,\,\,i = 0, 1,..., m - 1,\,\,p_m
= p_0.
$$
In other words, $\alpha$ has the tracial cyclic Rokhlin property.
\end{proof}

\vskip 3mm

\begin{thm}\label{IIT2}
{\rm
Let $A$ be a unital separable simple \CA\, with $\tr(A) = 0$ which has
a unique tracial state and
satisfies the Universal Coefficient Theorem.
Suppose that
$\alpha^r \in \Aut(A)$  is approximately inner
for some integer $r \geq 1$ and
that $\alpha^m$ is uniformly outer  for any integer $m \not= 0$.
Then $A \rtimes_{\alpha}{\mathbb Z}$ is a simple
AH-algebras with no dimension growth with real rank zero.
}
\end{thm}

\begin{proof}
Note that since $\alpha$ is outer,
$A \rtimes_{\alpha}{\mathbb Z}$ is simple by \cite{Ks0}.
{}From Theorem \ref{IT1} $\alpha$ has the tracial Rokhlin property.
Since $\alpha^r \in \Aut(A)$  is  approximately
inner for some integer $r \geq 1$, this implies that $\alpha$ has
the tracial cyclic Rokhlin property by Theorem \ref{IIT1}. So from
Theorem \ref{IT2} $\tr(A \rtimes_\alpha{\mathbb
Z})=0.$ Using the classification
theorem of \cite{Lnmsri} we  conclude that $A\rtimes_{\alpha}\Z$
is a simple AH-algebra with no dimension growth with real rank zero.
\end{proof}

\vskip 3mm

The following shows that the Kishimoto's conjecture that we
mentioned in the introduction is true at least for the case that
the simple $A{\mathbb T}$-algebra has a unique tracial state. In
Corollary \ref{IIC2}, we show that if one agrees that the
``sufficiently outer" means the automorphism has tracially Rokhlin
property then we do not need to assume that $A$ has a unique
tracial state.

\begin{Cor}\label{IIC1}
{\rm Let $A$ be a unital simple A${\mathbb T}$-algebra with a
unique trace and  real rank zero, and  let $\alpha\in \Aut(A)$
such that $\alpha$ is approximately inner.
 If
$\alpha^m$ is uniformly outer for any integer $m \not= 0$, or
$\alpha$ has tracial Rokhlin property, then $A
\rtimes_{\alpha}{\mathbb Z}$ is a unital simple A${\mathbb
T}$-algebra of real rank zero. }
\end{Cor}

\begin{proof}
From the following Pimsner-Voiculescu exact sequence\cite{PV},

$$
\begin{CD}
K_0 (A) @>{\id - \alpha^{-1}_*}>> K_0 (A) @>{\iota_*}>>
K_0 (A \rtimes_{\alpha}{\mathbb Z})\\
@AAA & &  @VVV\\
K_1 (A \rtimes_{\alpha}{\mathbb Z})
@<<{\iota_*}< K_1 (A) @<<{\id - \alpha^{-1}_*}< K_1 (A).
\end{CD}
$$

We sees  that $K_0(A \rtimes_{\alpha}{\mathbb Z})$ and
$K_1(A\times_{\alpha}{\mathbb Z})$ are torsion free. From
Theorem \ref{IT1}, Theorem \ref{IIT1}, and Theorem \ref{IT2}
we know that $\tr(A \rtimes_{\alpha}{\mathbb Z})=0$
and $A\rtimes_{\alpha}\Z$ satisfies the UCT. Therefore
$K_0(A\rtimes_{\alpha}\Z)$ is a weakly unperforated  Riesz group.
It follows from \cite{Ell2} that  there is a unital simple
A${\mathbb T}$-algebra $B$ with real rank zero 
which has the same ordered scaled
$K$-theory of $A \rtimes_{\alpha}{\mathbb Z}.$ It follows from
Theorem 5.1 of \cite{Lnmsri} that $A\cong B.$

\end{proof}

\begin{Cor}\label{IIC2}
{\rm Let $A$ be a unital simple A${\mathbb T}$-algebra (with real
rank zero)
 and $\alpha\in \Aut(A).$ Suppose that $\alpha$ is
approximately inner and $\alpha$ has the
tracial Rokhlin property. Then $A\rtimes_{\alpha}\Z$ is a unital
simple  A${\mathbb T}$-algebra  (with real rank zero.)}
\end{Cor}

\begin{proof}
It follows from Theorem \ref{IIT1} that $\alpha$ actually has the
tracial cyclic Rokhlin property. Then, by Theorem \ref{IT2},
$\tr(A\rtimes_{\alpha}\Z)\le 1$. As in the proof of \ref{IIC1},
$A\rtimes_{\alpha}\Z$ has torsion free $K$-theory. We then apply
the classification theorem in \cite{Lnmsri} (for real rank zero
case) or apply \cite{Lntai} (for real rank one case) to conclude that
$A\rtimes_{\alpha}\Z$ is a unital simple A\T-algebra  ( and with
real rank zero).
\end{proof}

\vskip 3mm

\begin{Remark}
{\rm Kishimoto in \cite{Ks1}, \cite{Ks2} and \cite{Ks3} proved
that if $A$ is a simple unital $A{\mathbb T}$-algebra of real rank
zero with a unique trace, and $\alpha \in \Aut(A)$ is an
approximately inner with the Rokhlin property, then  $A
\rtimes_\alpha{\mathbb Z}$ is also a simple unital $A{\mathbb
T}$-algebra under the assumption that both $K_0(A)$ and $K_1(A)$
are finitely generated with $K_1(A)\not= {\mathbb Z}$ and $\alpha
\in {\rm HInn}(A)$. Corollary \ref{IIC1} shows that the extra
conditions of $K_*(A)$ and $\alpha \in {\rm HInn}(A)$ are not
necessary. Corollary \ref{IIC2} shows that Kishimoto's conjecture
holds in general (without assuming that $A$ has the unique tracial
state) if the ``sufficiently outer" is replaced by the tracial
Rokhlin property. One should note that the tracial Rokhlin
property is weaker than the Rokhlin property used in
Kishimoto's work.(See Remark \ref{IR1}(iii).)
 Moreover, tracially cyclic Rokhlin property  is related to ``
 approximate Rokhlin" property in 4.2 of \cite{Ks1} which is also
weaker than the Rokhlin property used in Kishimoto's work. If one
allows the ``sufficiently outer" replaced by tracially cyclic
Rokhlin property, then $A\rtimes_{\alpha}\Z$ is always a unital
simple $A{\mathbb T}$-algebra without even assuming that
$\alpha^r$ is approximately inner but assuming
$A\rtimes_{\alpha}\Z$ has torsion free $K$-theory.}
\end{Remark}

\begin{thm}\label{IIITf}
Let $A$ be a unital separable simple \CA\, with $\tr(A)=0$ or
$\tr(A)=1$ and $\alpha\in \Aut(A)$ such that $\alpha^r$ is
approximately inner for some integer $r>0.$ Suppose that $\alpha$
has tracial Rokhlin property. Then $\tr(A\rtimes_{\alpha}\Z)=0,$ or
$\tr(A\rtimes_{\alpha}\Z)=1.$ Furthermore, if, in addition, $A$
satisfies the Universal Coefficient Theorem, then
$A\rtimes_{\alpha}\Z$ is a simple AH-algebra with no dimension
growth.
\end{thm}

\begin{proof}
The first part follows from Theorem \ref{IIT1} and Theorem
\ref{IT2}. For the last part, by \cite{Lntai}, $A$ is a simple
AH-algebra with no dimension growth. By the first part,
$\tr(A\rtimes_{\alpha}\Z)\le 1,$ it follows from \cite{Lntai} again that
$A\rtimes_{\alpha}\Z$ is also a simple AH-algebra with no
dimension growth.
\end{proof}

\begin{Remark}\label{RIIIF}
{\rm In Theorem \ref{IIT1}, we assume that $\tr(A)\le 1.$ 
In fact, we only need 
to assume that $A$ has the property (SP) and has the Fundamental Comparison
Property. Suppose that $A$ is a unital separable simple \CA\, 
with $\tr(A)=0$ and 
with a unique tracial state. Suppose also that 
$A\rtimes_{\alpha}\Z$ has a unique tracial state (unique ergodic).
Then by applying Theorem \ref{IIT1},
Theorem \ref{IT2} and Theorem \ref{IT1}, $\tr(A\rtimes_{\alpha}\Z)=0.$ 
On the other hand, in Corollary \ref{IIC1} and Corollary \ref{IIC2},
if we assume only that $\alpha^r$ is approximate inner 
(for $r>1$) and $\alpha$ has the tracial Rokhlin property, then 
$A\rtimes_{\alpha}\Z$ may not be an A${\mathbb T}$-algebra. 
This is because $A\rtimes_{\alpha}\Z$ may have torsion. 
However, it is a unital AH-algebra with no dimension growth 
by Theorem \ref{IIT2}.
But, in Corollary \ref{IIC1} and Corollary \ref{IIC2}, if we assume that
$\alpha^r$ is approximate inner for some integer $r$ and
$A\rtimes_{\alpha}\Z$ has torsion
free $K$-theory, then conclusion of both Corollary \ref{IIC1} and
Corollary \ref{IIC2} hold.
To allow torsion, related to the Kishimoto's conjecture, 
we proved (in  Theorem \ref{IIITf})  the following:
If $A$ is a unital simple AH-algebra with no dimension growth (with real rank
zero) and $\alpha\in \Aut(A)$ has the tracial Rokhlin property and 
$\alpha^r$ is approximate inner for some integer $r>0,$
then $A\rtimes_{\alpha}\Z$ is again a unital simple AH-algebra
with no dimension growth  (and with real rank zero).  
}

\end{Remark}

\section{Examples}

Let $G$ and $F$ be abelian groups. Recall that
$Pext(G,F)$ is the subgroup of those extensions
$$
0\to F\to E\to G\to 0
$$
so that each finitely generated subgroup of $G$ lifts.
If $A$ is a separable \CA\, which satisfies the  Universal Coefficient Theorem,
then, for any $\sigma$-unital \CA\, $B,$
$KL(A, B)=KK(A,B)/Pext(K_*(A), K_{*-1}(B)).$

\begin{Lem}\label{IVLKK}
{\rm
Let $A$ be a separable amenable \CA\, satisfying the UCT.
Suppose that $\alpha\in \Aut(A)$ such that
$(\alpha)_{*i}={\rm id}_{K_i(A)}.$ $i=0,1.$
Suppose that
$ext_{\Z}(K_{i-1}(A), K_i(A))/Pext(K_{i-1}(A), K_i(A))$ is finite.
Then there are  integers $r>0$ and $k>0$ such that
$$
[\alpha^{r+k}]=
[\alpha^k] \,\,\,{\rm in}\,\,\, KL(A, A).
$$
}
\end{Lem}

\begin{proof}
Consider $[\alpha^m]-[\alpha],$ for $m=1,2,....$
Since $(\alpha)_{*i}={\rm id}_{K_i(A)},$ $i=0,1,$
by the Universal Coefficient Theorem (\cite{RS}),
one computes that
$$
[\alpha^m]-[\alpha]\in ext_{\Z}(K_{*+1}(A), K_*(A)).
$$
Since $ext_{\Z}(K_{i-1}(A), K_i(A))/Pext(K_{i-1}(A), K_i(A))$ is
finite, there are positive integers $r$ and $k$ such that
$$
([\alpha^{r+k}]-[\alpha])=([\alpha^k]-[\alpha])\,\,\,{\rm
in}\,\,\, KL(A,A).
$$
It follows that
$$
[\alpha^{r+k}]=[\alpha^k]\,\,\,{\rm in}\,\,\,KL(A,A)
$$
\end{proof}

\begin{thm}\label{IVTainn}
{\rm Let $A$ be a unital separable simple \CA~ with $\tr(A) = 0$
satisfying the UCT. Suppose that $\alpha\in \Aut(A).$
 In any of
the following cases, $\alpha^r$ is approximately inner for some
integer $r>0.$ Consequently, if $\alpha^m$ is uniformly outer for
all $m\in \Z\backslash\{0\}$ (or $\alpha$ has tracial Rokhlin property),
$\alpha$ has tracial cyclic Rokhlin property and $\tr(A
\rtimes_\alpha{\mathbb Z}) = 0.$ In particular, $A\rtimes_{\alpha}
\Z$ is a simple AH-algebra with no dimension growth and real rank
zero.

\begin{itemize}
\item[$(1)$]
$K_0(A)=D,$ where $D$ is a countable dense subgroup of ${\mathbb R}$ and
$K_1(A)=\Z,$ or $K_1(A)=\{0\};$

\item[$(2)$] $K_0(A)=D,$ where $D$ is a finitely generated
countable dense subgroup of ${\mathbb R}$ and $K_1(A)=\Z$ or $K_1(A)$
is finite;

\item[$(3)$] $K_0(A)=D\oplus G,$
with $$
K_0(A)_+ = \{(r, x) \mid r \in D, r > 0, x \in G\} \cup \{(0,0)\},
$$
and $D$ is a dense subgroup of ${\mathbb R}$ such that for any
non-zero element $d\in D$ and any integer $n\ge 1,$ there is $e\in
D$ such that $me=d$ for some $m\ge n,$ where $G={\Z}$ or $G$ is
finite and $K_1(A)=\Z,$ or $K_1(A)=\{0\};$
\item[$(4)$]
$K_0(A) = \Q \oplus G$, where  $G={\Z}$ or $G$ is finite and $K
_1(A)=\Z,$ or $K_1(A)$ is a finite group.
\end{itemize}
}
\end{thm}

\begin{proof}
In all cases, it suffices to show that $\alpha^r$ is approximately inner
for some integer $r\ge 1.$

For (1), it is clear that $\alpha_{*0}={\rm id}_{K_0(A)}.$
If $K_1(A)=\Z,$
since $\alpha_{*1}$ is an isomorphism, $\alpha_{*1}(1)=\pm 1.$
Therefore $\alpha^2_{*i}={\rm id}_{K_i(A)}.$
Since $K_i(A)$ are torsion free, $[\alpha^2]=[{\rm id}_A]$ in $KL(A,A).$
It follows from Theorem 2.4 of \cite{Lnctaf} that $\alpha^2$ is approximately
inner.

For (2), as in (1), $\alpha_{*0}={\rm id}_{K_0(A)}.$
Also if $K_1(A)=\Z,$ then $\alpha^2_{*1}={\rm id}_{K_1(A)}.$
If $K_1(A)$ is finite, since $\alpha_{*1}$ is an isomorphism,
there exists $r_1\ge 1$ such that $\alpha^{r_1}_{*1}={\rm id}_{K_1(A)}.$
Let $\beta=\alpha^{2r_1}.$ Then $\beta_{*i}={\rm id}_{K_i(A)},$ $i=0,1.$
However, in this case,
$$
ext_{\Z}(D,K_1(A))=\{0\}\andeqn ext_{\Z}(K_1(A), K_0(A))\,\,\,{\rm
  is\,\,\, finite}.
$$
It follows from Lemma \ref{IVLKK} that
$[\beta^{m+k}]=[\beta^k]$ in $KL(A,A)$
for some integer $m,k\ge 1.$
By Theorem 2.3  of \cite{Lnctaf},
there exists a sequence of unitaries
such that
$$
\lim_{n\to\infty}{\rm ad}\, u_n\circ \beta^k(a)=\beta^{m+k}(a)\rforal a\in A.
$$
Since $\beta^k$ is an automorphism, it follows that $\beta^m$ is
approximately inner, or $\alpha^{m(2r_1)}$ is approximately inner.

For (3), as above, one has that $\alpha^2_{*1}={\rm id}_{K_1(A)}.$
The assumption on $D$ implies that there is no nonzero \hm\, from
$D$ to $\Z$ or a finite group. One then checks that there is an
integer $r_1\ge 1$ such that $\alpha_{*0}^{r_1}={\rm id}_{K_0(A)}.$ Put
$\beta = \alpha^{2r_1}.$ Then $\beta_{*i} = {\rm id}_{K_i(A)},$
$i=0, 1.$ To see that $\beta^m$ is approximately inner, we note
that $ext_{\Z}(D,K_1(A))=Pext(D,K_1(A))$ since $D$ is torsion
free. One then computes that
$$
ext_{\Z}(K_0(A), K_1(A))/Pext(K_0(A), K_1(A))= ext_
{\Z}(G,K_1(A))/Pext(K_0(A), K_1(A))
$$
which is finite, and $ext_{\Z}(K_1(A), K_0(A))=\{0\}.$ Thus one
can apply the same argument as in the case (2) by applying Lemma
\ref{IVLKK}.

For (4), as in the case (2) and (3), there is $r_1\ge 1$ such that
$\alpha^{r_1}_{*i}={\rm id}_{K_i(A)},$ $i=0,1.$ Moreover, since
$\Q$ is divisible,
$$
ext_{\Z}( K_1(A), \Q)=\{0\}.
$$
Because $K_1(A)=\Z,$ or $K_1(A)$ is finite,
$$
ext_{\Z}(K_1(A),K_0(A))=ext_{\Z}(K_1(A), G)\,\,\,{\rm
is\,\,\,also\,\,\, finite}.
$$
Since $\Q$ is torsion free, $ext_{\Z}(\Q, K_1(A))=Pext(\Q,
K_1(A)).$  One then computes that
$$
ext_{\Z}(K_0(A), K_1(A))/Pext(K_0(A), K_1(A))= ext_{\Z}(G,K_1(A))
/Pext(K_0(A), K_1(A))
$$
which is finite. Thus the argument in the case (3) applies.
\end{proof}

\vskip 5mm

\begin{Prop}\label{IIIP1}
{\rm Let $A$ be a simple unital \CA~ with $\tr(A) = 0$, and $B$ be
a simple unital \CA\, with $\tr(B) \leq 1$. Suppose that $\alpha
\in \Aut(A)$ has the tracial cyclic Rokhlin property. Then for any
$\beta \in \Aut(B)$ $\alpha \otimes \beta \in \Aut(A \otimes_{\rm
min}B)$ has the tracial cyclic Rokhlin property. }
\end{Prop}

\begin{proof}
It follows from \cite{HLX} $\tr(A \otimes_{\rm min }B) \leq 1$.
Hence $A \otimes_{\rm min }B$ has SP-property, stable rank one and
the Fundamental Comparison Property (\cite[Proposition 6.2 and
Theorems 6.9 and 6.11]{Lntr0}).

Let $F \subset A \otimes_{\rm min }B$ be a finite set, $n \in
{\mathbb N}$, and $\ep > 0$. Without loss of generality, we may
assume that there exist a finite set $F_A \subset A$ and $F_B
\subset B$ such that $F = F_A \otimes F_B$.

Since $\alpha$ has the tracial cyclic Rokhlin property, there exist
mutually orthogonal projections $e_0, e_1, \dots, e_n \in A$ such that

\begin{itemize}
\item[$(1)$]
$\| \alpha (e_j) - e_{j + 1} \| < \ep$ for $0 \leq j \leq n$,
where $e_{n + 1} = e_0$.
\item[$(2)$]
$\| e_j a - a e_j \| < \ep$ for $0 \leq j \leq n - 1$ and all $a \in F_A$.
\item[$(3)$]
$\tau(1 - \sum_{j = 0}^{n} e_j) < \ep$ for all tracial states $\tau$ on $A$.
\end{itemize}
(See Remark \ref{IR1} (ii).)

Set $f_i = e_i \otimes 1_B$ for $0 \leq i \leq n$.
Then $f_i$ are mutually orthogonal projections in
$A \otimes_{\rm min}B$
such that

\begin{itemize}
\item[$(1)$]
$\| (\alpha \otimes \beta)(f_j) - f_{j + 1} \| < \ep$ for $0 \leq j \leq n$,
where $f_{n + 1} = f_0$.
\item[$(2)$]
$\| f_j a - a f_j \| < \ep$ for $0 \leq j \leq n - 1$ and all $a \in F$.
\item[$(3)$]
$\tau(1 - \sum_{j = 0}^{n} f_j) < \ep$ for all tracial states $\tau$ on $A \otimes_{\rm min}B$.
\end{itemize}

This means that $\alpha \otimes \beta$ has the tracial cyclic Rokhlin property
by Remark \ref{IR1} (ii).
\end{proof}

\vskip 3mm

\begin{Cor}\label{IVCtr=1}
{\rm
Let $A$ be a separable simple amenable
unital \CA~ with $\tr(A) = 0$ which satisfies the UCT, and
   $B$ be a simple amenable unital \CA\, with $\tr(B) \leq 1$.
Suppose also that $A$ has a unique tracial state and
  $\alpha \in \Aut(A)$ such that
$\alpha^m$ is uniformly outer for all $m \not= 0$ and $\alpha^r$ is
approximately inner for some integer $r\ge 1.$
   Then for any $\beta \in \Aut(B),$
   $\alpha
\otimes \beta$ has the tracial cyclic Rokhlin property and
$\tr(D)\le 1,$ where $D = (A\otimes
B)\rtimes_{\alpha\otimes\beta}\Z.$ }
\end{Cor}

An unexpected consequence of the above corollary is that it
provides a new way to construct unital simple \CA s with tracial topological
rank one. All previous examples are inductive limit construction
(see \cite{G}). Since there is basically no restriction on $B$ and
$\beta,$ a  great number of those simple \CA s $D$ with $\tr(D) = 1$
can be obtained from Corollary \ref{IVCtr=1}. Since $\tr(B) = 1,$ one
certainly expects that most such $D$ has  $\tr(D)=1$ but not
$\tr(D)=0.$ To convince the reader that it is likely the case, we
compute the tracial rank in a very special case below. From its
construction, it should be clear how other example can be
constructed.

Denote $\Aff(A)$ the space of all affine continuous functions on
$\T(A)$. Given a projection $p \in M_n(A)$ for some integer $n
\geq 1$ we define $\rho_A(p)(\tau) = (\tau \otimes \Tr)(p)$ for
all $\tau \in \T(A)$, where $\Tr$ is the standard trace on
$M_n({\mathbb C})$. Then $\rho_A(p) \in \Aff(\T(A))$.

\begin{Ex}\label{IVE}
{\rm Let $A$ be a unital UHF-algebra with $K_0(A) = \Q$ and
$\alpha$ be in $\Aut(A)$ so that $\alpha^m$ is uniformly outer for
all $m \not= 0$ (or $\alpha$ is uniquely ergodic). Then $\alpha$ has
the tracial cyclic Rokhlin property by \cite[Lemma 4.3]{Ks1} and
Theorem \ref{IIT2}. Let $B$ be a unital simple $A{\mathbb
T}$-algebra for which $K_0(B) = \Q$ and $K_1(B) = \Z\oplus \Z$ and
$\Aff(\T(B)) = C_{\R}([0,1]).$
Existence of such simple
$A{\mathbb T}$-algebra was given by \cite{Tm}.
It follows from section 9
of  \cite{Tm} that there is $\beta\in \Aut(B)$ such that
$\beta_{*1}(x,y) = (-x,y)$ for $(x,y)\in \Z\oplus \Z,$
$\beta_{*0}= {\rm id}_{K_0(B)}$ and $\tau\circ \beta(b) = \tau(b)$
for all $b\in B$ and $\tau\in \T(B).$
It should be noted that
$\tau\circ \alpha(a) = \tau(a)$ for all $a\in \T(A).$ Put $\gamma =
\alpha\otimes \beta$ and $C = A\otimes B$ and $D =
C\rtimes_{\gamma}\Z.$

By the Kunneth formula one computes  that $C$ is a unital simple
($A{\mathbb T}$-algebra) with $K_0(C) = \Q$ and $K_1(C) = \Q\oplus \Q.$ One
computes that $\gamma_{*0}={ \rm id}_{K_0(C)}$ and
$\gamma_{*1}((x,y)) = (-x,y)$ for $(x,y)\in \Q\oplus \Q.$

It follows from Proposition \ref{IIIP1} that $\gamma$ has the
tracial cyclic Rokhlin property. Moreover $\tr(D) \le 1.$ To check
that $\tr(D) = 1,$ we first compute that, by Pimsner-Voiculescu's
exact sequence and by the divisibility of $\Q,$ we have $K_0(D) =
\Q\oplus \Q$ and $K_1(D) = \Q\oplus \Q.$ Consider  tracial states
with the form $t\otimes \tau,$ where $t\in \T(A)$ and $\tau \in
\T(B).$ Note all these tracial states are $\gamma$ invariant. Thus
they give tracial states on $D.$ Note that $\T(A)$ is a single
point. Thus we may identify $\T(B)$ with $\T(A)\otimes \T(B).$ Hence
$\Aff(\T(A)\otimes \T(B)) = C_{\R}([0,1]).$ Let $e_1=\rho_D((1,0)$
and $e_2=\rho_D((0,1)).$ Then
$$
\rho_D(K_0(D)) = \{ xe_1+ ye_2: x,y \in  \Q \}.
$$
We view $\T(B)\subset \T(D).$ Thus one has a surjective affine
\hm\, $\Lambda: \Aff(\T(D))\to \Aff(\T(B)) .$ It is easy to see
that $\Lambda\circ \rho_D(K_0(D))$ being rank two can not be dense
in $C_{\R}([0,1]).$ It follows from \cite[Theorem 6.9]{Bl} that
$D$ has real rank other than zero (actually one). It follows from
Theorem 7.1(c) of \cite{Lntr0} that $\tr(D)=1.$

If one insists to get non-zero torsion in $K$-theory,
one may start, for example, with $K_0(A)=\Q$ and
$K_1(A)=\Z/p\Z.$
}

\end{Ex}

\end{document}